\documentclass[12pt]{article}
\usepackage{amsmath, amsthm}\usepackage{enumerate}
\usepackage{amssymb}\usepackage{bbold}
\usepackage{color}
\usepackage[explicit]{titlesec}
\newtheorem{thm}{Theorem}
\newtheorem{prop}{Proposition}
\newtheorem{lem}{Lemma}

\newtheorem{cor}{Corollary}

\newtheorem{defi}{Definition}

\DeclareMathOperator{\ruc}{\xrightarrow[]{r}}

\renewcommand{\le}{\leqslant}
\renewcommand{\ge}{\geqslant}

\newcommand{\R}{\mathbb{R}}
\newcommand{\N}{\mathbb{N}}

\makeatletter
\renewcommand{\subsection}{\@startsection{subsection}{1}
{0pt}{3.25ex plus 1ex minus.2ex}{-1em}{\normalfont\normalsize\bf}}
\begin{document}
\title{Free uniformly complete vector lattices}
\date{\empty}
\maketitle
\author{\centering{Eduard Emelyanov and Svetlana Gorokhova\\ 
}}
\maketitle

\medskip
\medskip
\medskip


\medskip


%\medskip
%
%{\bf  Abstract:} {\normalsize We define a free uniformly complete 
%vector lattice $\text{\rm FUCVL}(A)$ over a non-empty set $A$ and give its representation 
%as a sublattice of the space $H(\mathbb{R}^A)$ of continuous in the product topology positively 
%homogeneous functions on $\R^A$.}

\section{Introduction and preliminaries}

\subsection{Free vector lattice over a set.}

We shall abbreviate a vector lattice (sublattice) by VL (resp., VSL).
Following the approach of de Pagter and Wickstead (cf. \cite[Definition 3.1]{PW}),
a {\em free vector lattice over} a non-empty set $A$ (briefly, $\text{\rm FVL}(A)$) 
is a pair $(F,i)$, where $F$ is a VL and $i:A\to F$ is a map such that, 
for any VL $E$ and for any map $q :A\to E$, there exists 
a unique Riesz homomorphism $T: F \to E$ satisfying 
$q=T \circ i$. If $(F,i)$ is a $\text{\rm FVL}(A)$ then $F$ is generated by $i(A)$ as a VL. 
The existence of a $\text{\rm FVL}(A)$ over $A$ is the long 
established fact going back to Birkhoff \cite{Birk}. A concrete representation 
of $\text{\rm FVL}(A)$ as a VL of real-valued functions was 
given by Weinberg \cite{Wein} and Baker \cite{Baker} 
(cf. also \cite{Blei}). Notice that $\text{\rm FVL}(A)$ is always Archimedean.

\smallskip
A {\em free Banach lattice} (briefly, $\text{\rm FBL}(A)$) {\em over 
a non-empty set} $A$ was introduced and investigated in \cite{PW} 
(see also \cite{Tr} for a simple alternative approach to a free 
lattice norm on $\text{\rm FVL}(A)$). In the present note, we define 
a {\em free uniformly complete vector lattice over} 
$A$ and give some of its representations. For further unexplained 
notation and terminology, we refer to \cite{Ku,PW}.

\subsection{Uniformly complete vector lattices.}
A net $x_\alpha$ in a VL $U$ {\em relatively uniformly converges} 
(briefly, {\em $\text{\rm r}$-converges}) 
to $x\in U$ if there exists $u \in U$ ({\em regulator of convergence}) 
such that, for each $k\in\N$, there exists $\alpha_k$ 
with $|x_\alpha-x|\le\frac{1}{k}|u|$ 
for all $\alpha\ge\alpha_k$ 
(in this case, we write $x_\alpha\ruc x(u)$). 
A net $x_\alpha$ in a VL $U$ is {\em $\text{\rm r}$-Cauchy} with a regulator $u\in U$
(briefly, $\text{\rm r}(u)$-{\em Cauchy}) if, for each $k\in\N$, there exists $\alpha_k$
with $|x_{\alpha'}-x_{\alpha''}|\le\frac{1}{k}|u|$ for all $\alpha',\alpha''\ge\alpha_k$.

\smallskip
A VSL $V$ of a VL $U$ is called \text{\rm r}-{\em complete in} $U$ if
$U$ is Archimedean and, for every $u\in U$ and every $\text{\rm r}(u)$-Cauchy
net $x_\alpha$ in $V$, there exist $x\in V$ and $w\in U$ 
with $x_\alpha\ruc x(w)$ (such an $x$ is unique since $V$ is Archimedean). 
For example, VSL $c_{00}$ of $c_0$ is not \text{\rm r}-complete in $c_0$ as 
$c_{00}\ni\sum_{i=1}^{n}\frac{1}{k^2}e_i\ruc
\sum_{i=1}^\infty\frac{1}{k^2}e_i\notin c_{00}$
with regulator $\sum_{i=1}^\infty\frac{1}{k}e_i\in c_0$, 
yet $c_{00}$ is \text{\rm r}-complete in itself.

\smallskip
We include the following, certainly well known fact, for which we did not 
found an appropriate reference.

\begin{lem}[]\label{L1-0}
Let $x_\alpha$ be an $\text{\rm r}(u)$-Cauchy
net in an Archi\-me\-dean \text{\rm VL} $U$, where $u\in U$.
If $x_\alpha\ruc x(w)$ with $x,w\in U$, then $x_\alpha\ruc x(u)$.
\end{lem}

\begin{proof}
For each $l\in\N$ take an $\alpha(l)$ such that $|x_\alpha-x|\le\frac{1}{l}|w|$ 
when $\alpha\ge\alpha(l)$. 
Since $x_\alpha$ is $\text{\rm r}(u)$-Cauchy,
for each $k\in\N$, there exists $\alpha_k$ with 
$|x_{\alpha'}-x_{\alpha''}|\le\frac{1}{k}|u|$
for $\alpha',\alpha''\ge\alpha_k$. Now, for every $l,k\in\N$, 
we pick $\alpha(k,l)\ge\alpha_k,\alpha(l)$. 
Then, for each $\alpha\ge\alpha_k$, 
$$
   |x_\alpha-x|\le|x_\alpha-x_{\alpha(k,l)}|+
   |x_{\alpha(k,l)}-x|\le\frac{1}{k}|u|+\frac{1}{l}|w|.
   \eqno(1)
$$
Since $l$ is arbitrary in (1) and $U$ is Archimedean then
$|x_\alpha-x|\le\frac{1}{k}|u|$ for $\alpha\ge\alpha_k$,
and hence $x_\alpha\ruc x(u)$.
\end{proof}

\begin{cor}[]\label{same reg}
Let $V$ be a \text{\rm VSL} of an Archimedean \text{\rm VL} $U$.
The following conditions are equivalent.
\begin{enumerate}[$i)$]
\item
$V$ is $\text{\rm r}$-complete in $U$.
\item
For every $u\in U$ and every $\text{\rm r}(u)$-Cauchy net $x_\alpha$ in $V$,
there exists $x\in V$ with $x_\alpha\ruc x(u)$.
\end{enumerate}
\end{cor}

\smallskip
As the \text{\rm r}-convergence is sequential
(cf., e.g. \cite[Proposition 2.4]{AEG}), 
we can restrict ourselves to \text{\rm r}-convergent/\text{\rm r}-Cauchy sequences. 
In particular, a VSL $V$ of an Archimedean VL $U$ is \text{\rm r}-complete in $U$
iff, for every $u\in U$ and every $\text{\rm r}(u)$-Cauchy sequence $x_n$ in $V$,
there exists $x\in V$ such that $x_n\ruc x(u)$. 

\smallskip
An \text{\rm r}-complete in itself VL will be abbreviated by \text{\rm UCVL}.
For example, Dedekind complete VLs and Banach lattices are UCVLs.
Note that each UCVL is Archimedean by definition.

\begin{lem}[]\label{princ-ideal}
Let $U$ be a \text{\rm VL}. The following conditions are equivalent.
\begin{enumerate}[$i)$]
\item
$U$ is a \text{\rm UCVL}.
\item
For every $y\in U$, the principal order ideal $U_y=\bigcup_{n\in\mathbb{N}}[-n|y|,n|y|]$ 
is complete under the norm $\|x\|_y:=\inf\{\lambda\in\mathbb{R}:|x|\le\lambda|y|\}$.
\item
$U_y$ is \text{\rm UCVL} for all $y\in U$.
\end{enumerate}
\end{lem}

\begin{proof}
$i)\Longrightarrow ii)$\
Let $y\in U$. As $U$ is Archimedean, the map $x\to\|x\|_y$ is a norm on $U_y$.
Since $\|x_n-x\|_y\to 0\Longleftrightarrow x_n\xrightarrow[]{r}x(|y|)$ 
for each sequence $x_n$ in $U$ and since $U$ is \text{\rm r}-complete, then $U_y$ 
is $\|\cdot\|_y$-complete.

\smallskip
$ii)\Longrightarrow i)$\
Since the map $x\to\|x\|_y$ is a norm on $U_y$ for every $y\in U$, then $U$ is Archimedean.
Let $u\in U$ and let $x_n$ be an $\text{\rm r}(u)$-Cauchy sequence in $U$.
Then, there exists $n_1$ satisfying $|x_k-x_m|\le|u|$ for all $k,m\ge n_1$.
Therefore, $|x_k-x_m|\le w:=|u|+\sum_{i=1}^{n_1}|x_i|\in U$ for all $k,m\in\mathbb{N}$,
and hence the sequence $x_n$ lies in $U_w$. Since $x_n$ is  $\text{\rm r}(u)$-Cauchy 
then it is  $\text{\rm r}(w)$-Cauchy,
and hence it is $\|\cdot\|_w$-Cauchy in $U_w$. Since $U_w$ is $\|\cdot\|_w$-complete,
there exists $x\in U_w$ with $\|x_n-x\|_w\to 0$, and hence $x_n\ruc x(w)$.
Then $x_n\ruc x(u)$ by Lemma \ref{L1-0}. 

\smallskip
$ii)\Longleftrightarrow iii)$\
It is trivial.
\end{proof}

\begin{defi}[]\label{d00}
A \text{\rm UCVL} $F$ is an \text{\rm r}-{\em completion} 
of a \text{\rm VL} $E$ if there exists a Riesz embedding $i: E\to F$ such that, 
for each \text{\rm UCVL} $G$ and each Riesz homomorphism $T: E\to G$ 
there is unique Riesz homomorphism $S: F\to G$ satisfying $T=S\circ i$.
\end{defi}

Since every UCVL is Archimedean, each VL possessesing an \text{\rm r}-completion
is Archimedean. 

\smallskip
If an $\text{\rm r}$-completion $F$ of a VL $E$
exists, it must be unique up to a Riesz isomorphism. Indeed,
let a UCVL $F_1$ be another $\text{\rm r}$-completion of $E$. 
Then, there are two 
Riesz embeddings $i: E\to F$ and $i_1: E\to F_1$ with two 
Riesz homomorphisms $S: F\to F_1$ and $S_1: F_1\to F$ 
satisfying $i_1=S\circ i$ and $i=S_1\circ i_1$
and hence $i=S_1\circ S\circ i$. 
Since a Riesz homomorphism $H: F\to F$ satisfying $i=H\circ i$ is unique then $S_1\circ S=I_F$,
where $I_F$ is the identity map on $F$.
Similarly, $S\circ S_1=I_{F_1}$. Thus, VLs $F$ and $F_1$ are Riesz isomorphic.

\smallskip
It is long known (e.g., see Veksler's works \cite{Veks1,Veks2}) that the intersection 
of all \text{\rm r}-complete VSLs containing $E$ in the Dedekind completion 
$E^\delta$ of an Archimede\-an VL $E$ is an $\text{\rm r}$-completion of $E$. 
For convenience of the reader, we include details of the construction of $\text{\rm r}$-completion 
in slightly more general setting. 

\subsection{$\text{\rm r}$-Completion of an Archimedean vector lattice.}
Throughout this subsection, $E$ is a \text{\rm VSL} of a \text{\rm UCVL} $U$.

\begin{lem}\label{Rem1}
The intersection $F$ of all \text{\rm UCVL}s in $U$ containing $E$ is a \text{\rm UCVL}.
\end{lem}

\begin{proof}
Suppose $x_n$ is an $\text{\rm r}(x_1)$-Cauchy sequence in $F$. 
Take a UCVL $V$ in $U$ containing $E$. Then $x_n\ruc x(x_1)$ for some $x\in V$ 
by Corollary \ref{same reg}. Since $V$ is arbitrary UCVL in $U$ containing $E$, 
and since $\text{\rm r}$-limit $x$ of the sequence $x_n$ is unique, $x$ belongs 
to all UCVLs in $U$ containing $E$, and hence $x\in F$.
\end{proof}

Define a transfinite sequence $(F_\beta)_{\beta\in{\text{\rm Ord}}}$ of VSLs of $U$ by:
\begin{enumerate}[-]
\item 
$F_1:=E$;
\item
$F_{\beta+1}:=\{x\in U: x_n \ruc x(x_1),\ 
\text{\rm for a sequence}\ x_n\ \text{\rm in}\ F_{\beta}\}$;
\item
$F_\beta:=\bigcup_{\gamma\in{\text{\rm Ord}};\ \gamma<\beta}F_{\gamma}$
for a limit ordinal $\beta$.
\end{enumerate}

Note that $F_{\beta_1}$ is a VSL of $F_{\beta_2}$ if $\beta_1\le\beta_2$.

\begin{lem}[]\label{L1}
The intersection $F$ of all \text{\rm UCVL}s of $U$ containing $E$ satisfies
$F=\bigcup_{\gamma\in{\text{\rm Ord;$\gamma<\omega_1$}}}F_{\gamma}$.
\end{lem}

\begin{proof}
Take an $\text{\rm r}(x_1)$-Cauchy sequence $x_n$ in 
$G:=\bigcup_{\gamma\in{\text{\rm Ord;$\gamma<\omega_1$}}}F_{\gamma}$. 
Since $U$ is a UCVL, $x_n\ruc x(x_1)$ for some $x\in U$ by Lemma \ref{L1-0}. 
Pick an ordinal $\beta<\omega_1$ with $x_n\in F_\beta$ for all $n$.
Then $x\in F_{\beta+1}\subseteq G$. 
Therefore, $G$ is a UCVL in $U$ containing $E$, and hence $F\subseteq G$ due to
Lemma \ref{Rem1}.

\smallskip
Let $V$ be a UCVL of $U$ containing $E$, so $F_1=E\subseteq V$. 
Let $\gamma$ be an ordinal such that $F_{\beta}\subseteq V$ for all $\beta<\gamma$.
If $\gamma$ is a limit ordinal then
$F_\gamma=\bigcup_{\beta<\gamma}F_{\beta}\subseteq V$.
If $\gamma=\beta+1$ then, for each $x\in F_\gamma$,
there exists a sequence $x_n$ in $F_\beta\subseteq V$
with $x_n \ruc x(x_1)$. As $V$ is UCVL then $x\in V$,
and since $x\in F_\gamma$ is arbitrary, $F_\gamma\subseteq V$.
In both cases, $G\subseteq\bigcup_{\gamma\in{\text{\rm Ord}}}F_{\gamma}\subseteq V$.
Since $V$ is arbitrary UCVL of $U$ containing $E$, then $G\subseteq F$
by Lemma \ref{Rem1}. In view of the already obtained inclusion
$F\subseteq G$, we conclude $F=G$.
\end{proof}

\begin{prop}[]\label{L2}
Let $E$ be a \text{\rm VSL} of a \text{\rm UCVL} $U$. 
Then the intersection $F$ of all \text{\rm UCVL}s of $U$ containing $E$
is an $\text{\rm r}$-completion of $E$.
\end{prop}

\begin{proof}
By Lemma \ref{L1}, $F$ is a UCVL that coincides with 
$\bigcup_{\gamma\in{\text{\rm Ord;$\gamma<\omega_1$}}}F_{\gamma}$.
In order to show that $F$ is an $\text{\rm r}$-completion of $E$,
let $G$ be a UCVL and let $T: E\to G$ be a Riesz 
homomorphism. It is enough to show that $T$ has
a unique extension to a Riesz homomorphism $S: F\to G$. 

\smallskip
First we prove that such an extension $S$,
if exists, is unique. Indeed assume $S_1: F\to G$ is another such an extension.
Clearly, $S$ agrees with $S_1$ on $F_1=E$. 
Assume $S$ agrees with $S_1$ on $F_{\beta}$ for all $\beta<\gamma$.
Let $x\in F_{\gamma}$.
If $\gamma$ is a limit ordinal then
$x\in F_{\gamma}=\bigcup_{\beta<\gamma}F_{\beta}$.
Then $x\in F_{\beta}$ for some $\beta<\gamma$, and hence
$Sx=S_1x$ by the assumption.
If $\gamma=\beta+1$ then there exists a sequence 
$x_n$ in $F_\beta$ with $x_n \ruc x(x_1)$. In other words, for each $k\in\mathbb{N}$, 
there exists $n_k$ with $|x_n-x|\le\frac{1}{k}|x_1|$ for all $n\ge n_k$.
Then 
$$
   |Sx_n-Sx|+|S_1x_n-S_1x|\le\frac{1}{k}(S|x_1|+S_1|x_1|)\ \ \ \ \ (\forall n\ge n_k).
   \eqno(2)
$$
Because the sequence $x_n$ is contained in $F_\beta$ then $Sx_n=S_1x_n$ 
for all $n\in\mathbb{N}$. Now, it follows from (2) that
$$
   Sx_n\ruc Sx(S|x_1|+S_1|x_1|) \ \ \ \text{\rm and} \ \ \ 
   S_1x_n\ruc S_1x(S|x_1|+S_1|x_1|)
$$ 
in $G$, and since $G$ is Archimedean, then $Sx=S_1x$ as desired.

\smallskip
In order to show the existence of the required extension $S: F\to G$,
we use Lemma \ref{L1} and apply the transfinite induction by $\gamma$.
Suppose the extensions $S_\alpha: F_\alpha\to G$ exist already
for all $\alpha<\gamma$. By the arguments as above,
they are unique for all $\alpha<\gamma$,
and hence $S_\alpha$ agree with $S_\beta$ on $F_{\min{(\alpha,\beta)}}$.
Let $x\in F_\gamma$.

\smallskip
(A) If $\gamma$ is a limit ordinal then $x\in F_\alpha$ 
for some $\alpha<\gamma$. We set $S_\gamma x:=S_\alpha x$.

\smallskip
(B) If $\gamma=\beta+1$, there exists a sequence 
$x_n$ in $F_\beta$ with $x_n \ruc x(x_1)$. Then, for each $k\in\mathbb{N}$, 
there exists $n_k$ with $|x_n-x|\le\frac{1}{k}|x_1|$ for all $n\ge n_k$, and hence 
$$
   |x_n-x_m|\le|x_n-x|+|x_m-x|\le\frac{2}{k}|x_1| \ \ \ \ (\forall n,m\ge n_k).
   \eqno(3)
$$
It follows from (3) that $|S_\beta x_n-S_\beta x_m|\le\frac{2}{k}S_\beta|x_1|\in G$ 
for all $n,m\ge n_k$.
So, the sequence $S_\beta x_n$ is $\text{\rm r}$-Cauchy in $G$, and hence there exists 
unique $y\in G$ with $S_\beta x_n\ruc y$ in $G$. 
It is easy to see that $y$ does not depend on the choice of the sequence 
$x_n$ in $F_\beta$ such that $x_n \ruc x(x_1)$. We set $S_\gamma x:=y$.

\smallskip
Finally, let $x\in F$. Then $x\in F_\gamma$ for some $\gamma$. 
We set $Sx:=S_\gamma x$. Since linear and lattice operations in VLs
are $\text{\rm r}$-continuous, the constructed extension $S: F\to G$ 
is a Riesz homomorphism. 
\end{proof}

\smallskip
Proposition \ref{L2} combined with Lemma \ref{princ-ideal} and with the
theorem of Kakutani and M. Krein -- S. Krein gives that each principal order ideal
of an $\text{\rm r}$-completion of an Archimedean VL is Riesz isomorphic
to $C(K)$ for some compact Hausdorff space $K$ (cf., \cite[Theorem 2]{Veks1}).

\subsection{Free \text{\rm r}-complete vector lattice over a set.}
Although the next definition does not, as such, 
appear in the literature, we think it belongs to folklore.

\begin{defi}[]\label{d1}
If $A$ is a non-empty set, then a {\em free UCVL over} $A$ 
{\em (}briefly, $\text{\rm FUCVL}(A)${\em )} is 
a pair $(F,f)$, where $F$ is a \text{\rm UCVL} and $f:A\to F$ is an embedding 
such that, for any \text{\rm UCVL} $E$ and for any embedding $q:A\to E$, 
there exists a unique Riesz homomorphism $T:F\to E$ satisfying $q=T\circ f$. 
\end{defi}

If exists, $\text{\rm FUCVL}(A)$ is an initial object in 
the category, whose objects are pairs $(E,q)$ with a UCVL $E$ and an embedding $q :A\to E$, 
and whose morphisms $q_2^1:(E_1, q_1)\to(E_2, q_2)$ are Riesz homomorphisms 
from $E_1$ to $E_2$ satisfying $q_2=q_2^1\circ q_1$. 
Thus, $\text{\rm FUCVL}(A)$ is defined similarly to $\text{\rm FVL}(A)$ in a proper subcategory. 
Routine arguments show that a free UCVL $(E,q)$ over $A$, if exists,
is unique up to a Riesz isomorphism.

\begin{thm}\label{maintt1}
Let $A$ be a non-empty set, $\text{\rm FVL}(A)$ be a \text{\rm VSL} 
of a \text{\rm UCVL} $U$, and $A\xrightarrow[]{i}\text{\rm FVL}(A)$ be 
the assigned embedding. Then$:$
\begin{enumerate}[$(i)$]
\item
the intersection $F$ of all \text{\rm UCVL}s of $U$ containing $\text{\rm FVL}(A)$  
together with the embedding $i$ is a $\text{\rm FUCVL}(A)$$;$
\item
$F=\bigcup\limits_{y\in\text{\rm FVL}(A)}U_y=
\bigcup\limits_{B\ \text{\rm is a nonempty finite 
subset of}\ A}U_{\sum\limits_{b\in B}|i(b)|}$.
\end{enumerate}
\end{thm}

\begin{proof}
$(i)$\
Note that $F$ is an $\text{\rm r}$-completion 
of $\text{\rm FVL}(A)$ due to Proposition \ref{L2}.

\smallskip
Let $E$ be a UCVL and let $q: A \to E$ be an embedding. 
By the definition of $\text{\rm FVL}(A)$, there exists 
a unique Riesz homomorphism $T_1:\text{\rm FVL}(A)\to E$
with $q = T_1\circ i$. Since $F$ is an $\text{\rm r}$-completion 
of $\text{\rm FVL}(A)$, there is 
unique Riesz homomorphism $T:F\to E$ satisfying $T_1 = T\circ J$,
where $J$ is the inclusion of $\text{\rm FVL}(A)$ into $F$, and hence
$q = T_1\circ i=T\circ J\circ i=T\circ i$.
By Definition \ref{d1}, the pair $(F,i)$ is a $\text{\rm FUCVL}(A)$.

\smallskip
$(ii)$\
It follows directly from Lemma \ref{princ-ideal}.
\end{proof}

Let $A$ be a non-empty set. Recall that $\text{\rm FVL}(A)$ is the VSL of 
UCVL $\R^{\R^A}$ generated by the evaluation functionals $\delta_a$ 
on $\R^A$, $\delta_a(\xi)=\xi(a)$ \cite[Theorem 2.4]{Baker}. The next 
representation of $\text{\rm FUCVL}(A)$ follows from Theorem \ref{maintt1}.

\begin{cor}\label{cortt1}
Let $A$ be a non-empty set and let $E$ be the \text{\rm VSL} of a \text{\rm UCVL} 
$\R^{\R^A}$ generated by evaluation functionals $\delta_a$, $a\in A$
with the embedding $A\xrightarrow[]{i}\R^{\R^A}$, $i(a):=\delta_a$. 
Then, the intersection $F$ of all \text{\rm UCVL}s of $\R^{\R^A}$ containing 
$i(A)$ together with embedding $i$ is a $\text{\rm FUCVL}(A)$.
\end{cor}

\subsection{Restrictions and projections.}
Let $A$ be a nonempty set. Following the tradition, for $\emptyset\ne B\subseteq A$, we identify 
$\R^{\R^B}$ with a VSL of  $\R^{\R^A}$
by assigning $\xi \in\R^{\R^B}$ with $\hat{\xi} \in\R^{\R^A}$ 
such that $\hat{\xi} (f) = \xi(f|_B)$ (see \cite{PW}).
By \cite[Proposition~3.5(2)]{PW}, there exists a unique Riesz 
homomorphism projection $P_B$ of $\text{\rm FVL}(A)$  onto $\text{\rm FVL}(B)$ satisfying
$$
   P_B(\delta_a)=
   \begin{cases}
   \delta_a  &  \text{if}  \ \  a\in B \\
   0  & \text{if}  \ \  a\in A\setminus B \, .
   \end{cases}
$$
In particular, 
$
   \text{\rm FVL}(A)= \bigcup\limits_{\emptyset\ne B\in{\cal P}_{fin}(A)}
   \text{\rm FVL}(B)\subseteq\R^{\R^A} 
$,
where ${\cal P}_{fin}(A)$ is the set of all finite subsets of $A$ (see \cite[Proposition 3.7]{PW}). 

\smallskip
Denote by $H(\R^A)$ (by $H(\Delta_A)$) the space of all positively 
homogeneous real-valued functions on $\R^A$ (on $\Delta_A:=[-1,1]^A$) 
which are continuous in the product topology of $\R^A$ (of $\Delta_A$). 
Clearly, $H(\Delta_A)$ is a VSL of $C(\Delta_A)$, closed in $\|.\|_\infty$. 
Thus, $H(\Delta_A)$ is a Banach lattice, and hence is a UCVL. 
It is well known that $\text{\rm FVL}(A)$ may be identified with a VSL of $H(\R^A)$ 
and hence with a VSL of $H(\Delta_A)$ in view of \cite[Lemma~5.1]{PW}. 
By \cite[Corollary~5.7]{PW}, $\text{\rm FBL}(A)$ is embedded into $H(\Delta_A)$ as 
a VSL $J(\text{\rm FBL}(A))$ and even as an order ideal, as mentioned in \cite{PW}. 
Furthermore,  $J(\text{\rm FBL}(A))=H(\Delta_A)$ iff $A$ is finite and, in this case,  
$\text{\rm FBL}(A)$ is isomorphic to $H(\Delta_A)$ under 
the supremum norm \cite[Theorem~8.2]{PW}.

\begin{thm}\label{l2}
Let $B$ be a non-empty finite set. Then $\text{\rm FUCVL}(B)$ is Riesz isomorphic 
to $\text{\rm FBL}(B)$, $H(\Delta_B)$, and to $H(\R^B)$. 
\end{thm}

\begin{proof}
As $\text{\rm FVL}(B)$ and $\text{\rm FBL}(B)$ possess (the same) strong 
order unit (cf. \cite[Proposition~5.3 and Theorem~8.2]{PW}) and $\text{\rm FBL}(B)$ 
is Riesz isomorphic to the completion of $\text{\rm FVL}(B)$ under its strong unit norm, 
$\text{\rm FBL}(B)$ is the $\text{\rm r}$-closure of $\text{\rm FVL}(B)$,
and hence $\text{\rm FUCVL}(B)$ is Riesz isomorphic to $\text{\rm FBL}(B)$.
The fact that $\text{\rm FBL}(B)$ 
and $H(\Delta_B)$ are Riesz isomorphic is the result of \cite[Proposition~5.3]{PW}. 
As $B$ is finite, the restriction map $R: H(\R^B) \to H(\Delta_B)$ is 
surjective and hence is a Riesz isomorphism (e.g. by \cite[Lemma~5.1]{PW}).
\end{proof}

For the convenience, we include the proof of the following certainly well known proposition.

\begin{prop}\label{p2}
Let $A\ne\emptyset$, and let $x\in\text{\rm FBL}(A)$. 
Then  there exists a sequence $x_n$ in $\text{\rm FVL}(A)$ 
that $\text{\rm r}$-converges to $x$ with a regulator $u\in\text{\rm FBL}(A)$.
\end{prop}

\begin{proof}
Since  $x\in\text{\rm FBL}(A)$, then there exists a sequence $x_n$ in $\text{\rm FVL}(A)$ 
satisfying $\|x_n - x\|_F \le \frac{1}{n2^n}$ for all $n\in \N$, where $\|.\|_F$ is 
the extension to $\text{\rm FBL}(A)$ of a free lattice norm on $\text{\rm FVL}(A)$. 
Let $u=\|.\|_F\text{-}\sum_{n=1}^{\infty}n |x_n -x|\in\text{\rm FBL}(A)$. 
Then $|x_n -x| \le\frac{1}{n}u$  for all $n\in \N$
and hence $x_n \ruc x(u)$.
\end{proof}

Since any Banach lattice is UCVL, it follows from Proposition \ref{p2} 
that $\text{\rm FBL}(A)$ contains $\text{\rm r}$-completion of $\text{\rm FVL}(A)$. 
We show that $\text{\rm FUCVL}(A)$ is a proper VSL 
of $\text{\rm FBL}(A)$ unless $A$ is finite. 

\begin{prop}\label{p3}
Let $A\ne\emptyset$. 
If a sequence $g_n \in H(\R^A)$  $\text{\rm r}$-converges with a regulator  
$u \in H(\R^A)$ to some  $g \in \R^{\R^A}$ then $g \in H(\R^A)$. 
\end{prop}

\begin{proof}
Without lost of generality, we may assume that $|g_n - g| \le \frac{1}{n} u$ for all $n\in \N$. 
Clearly, $g$ is positively homogeneous. It remains to show the continuity of $g$. 
Let $\eta_\alpha \to \eta$  in the product topology $\tau$ of $\R^A$. 
Then, for some $M$, there is $\alpha_0$ with $u(\eta_\alpha) \le M$  
for all $\alpha \ge \alpha_0$. Hence
$$
  |g(\eta_\alpha) - g(\eta)| \le |g(\eta_\alpha) -g_n(\eta_\alpha) |+ 
  |g_n(\eta_\alpha) -g_n(\eta) |+ |g_n(\eta) - g(\eta)|\le 
$$
$$
   |g_n(\eta_\alpha) -g_n(\eta) |+\frac{2M}{n}  
  \quad\quad  (\forall n\in \N) \quad (\forall \alpha \ge \alpha_0).
$$
Since $g_n$ is $\tau$-continuous, 
there exists $\alpha_n$ with $|g_n(\eta_\alpha) -g_n(\eta) |\le \frac{M}{n}$ 
for $\alpha \ge \alpha_n$. 
It follows that $|g(\eta_\alpha) - g(\eta)| \le \frac{3M}{n}$, 
if $\alpha \ge \alpha_0, \alpha_n$. Thus $g$ is $\tau$-continuous and hence $g \in H(\R^A)$.
\end{proof}

\begin{cor}\label{new}
$H(\R^A)$ is an $\text{\rm r}$-closed \text{\rm VSL} of $\R^{\R^A}$ 
with regulators from $H(\R^A)$.
\end{cor}

Considering $\text{\rm FVL}(A)$ as a VSL of $H(\Delta_A)$, 
we obtain the following result.

\begin{cor}\label{duble}
Let $R: H(\R^A)\to H(\Delta_A)$ be the restriction map.
If a sequence $g_n$ in $\text{\rm FVL}(A)$  $\text{\rm r}$-converges 
to some $g\in H(\Delta_A)$ with a regulator $u\in R(H(\R^A))$ then 
$g \in R(H(\R^A))$,  
\end{cor}

It is worth to compare Corollary \ref{duble} with  \cite[Example~5.2]{PW}. 

\subsection{Applications.}
The following theorem is a refinement of Theorem~\ref{l2}.

\begin{thm}\label{T1-1}
Let $A$ be a non-empty set. The intersection $F$ of all \text{\rm r}-closed 
\text{\rm VSL}s of $H(\R^A)$ containing $\text{\rm FVL}(A)$ is a $\text{\rm FUCVL}(A)$.
\end{thm}

\begin{proof}
By Corollary~\ref{new}, $H(\R^A)$ is a VSL of $\R^{\R^A}$
\text{\rm r}-closed under the regulators of $H(\R^A)$.
Since $\R^{\R^A}$ is a UCVL then $H(\R^A)$ is a UCVL.
As $\text{\rm FVL}(A)$ is a VSL of $H(\R^A)$, then
$\text{\rm FUCVL}(A)$ is an \text{\rm r}-closed in $H(\R^A)$ VSL of $H(\R^A)$.
\end{proof}

\begin{cor}\label{CorB}
Let $A$ be a non-empty set. Then 
$$
  \bigcup\limits_{\emptyset\ne B\in{\cal P}_{fin}(A)}\text{\rm FBL}(B)\subseteq 
  \text{\rm FUCVL}(A)\subseteq\text{\rm FBL}(A).
  \eqno(4)
$$
Furthermore, both inclusions are proper unless $A$ is finite. 
\end{cor}

\begin{proof}
Formula (4) follows from Lemma \ref{princ-ideal},
\cite[Theorem 8.2]{PW}, and Theorem \ref{l2} 
since every \text{\rm FBL} is an $\text{\rm UCVL}$. 
Let $A$ be infinite. The inclusion $\text{\rm FUCVL}(A)\subseteq\text{\rm FBL}(A)$ 
is proper since e.g. 
$g: = \sum_{k=1}^{\infty}2^{-k}\delta_{a_k}\in\text{\rm FBL}(A)$ for distinct $a_k \in A$
yet $g \not\in H(\R^A)$ (cf. \cite[Example~5.2]{PW}).
The inclusion $\bigcup\limits_{\emptyset\ne B\in{\cal P}_{fin}(A)}
\text{\rm FBL}(B)\subseteq\text{\rm FUCVL(A)}$
is proper since $z: = \sum\limits_{k=1}^\infty 2^{-k}|\delta_{a_k}|\wedge|\delta_{a_1}| 
\in\text{\rm FUCVL}(A)$ for distinct $a_k\in A$ yet $z\not\in\text{\rm FBL}(B)$ 
for each finite nonempty $B\subseteq A$.
\end{proof}

\end{document}